\documentclass{article}
\usepackage{amsmath,amsthm,amssymb,array}
\usepackage[all]{xy}

\theoremstyle{plain}
\newtheorem{theorem}{Theorem}
\newtheorem{lemma}{Lemma}

\theoremstyle{definition}
\newtheorem{example}{Example}

\newcommand\Aut{{\mathsf{Aut}}}
\newcommand\Stab{{\mathsf{Stab}}}
\newcommand\Sym{{\mathsf{S}}}
\newcommand\tree{{\mathcal T}}
\newcommand\G{{\mathcal G}}
\newcommand\GG{{\bar{\mathcal G}}}

\title{Construction of elements in the closure of Grigorchuk group}
\author{Goulnara Arzhantseva and Zoran \v{S}uni\'c}
\date{Universit\'e de Gen\`eve, Section de Math\'ematiques, 2-4, rue du
Li\`evre, Case postale 64, CH-1211 Gen\'eve 4, Suisse
\vspace{2mm}\\Department of Mathematics, Texas A\&M University,
College Station, TX 77843-3368, USA}

\begin{document}
\maketitle

%%%%%%%%%%%%%%%%%%%%%%%%%%%%%%%%%%%%%%%%%%%%%%%%%%%%%

The group known as Grigorchuk group (also the first Grigorchuk
group) was introduced in~\cite{grigorchuk:burnside}. More
information on this remarkable group can be found
in~\cite{grigorchuk:unsolved} and in~\cite{harpe:ggt}. Here we only
introduce as much as is necessary to describe the elements in the
closure $\GG$ of Grigorchuk group $\G$ in the pro-finite group
$\Aut(\tree)$ of binary rooted tree automorphisms. In fact, we
describe constraints that need to be satisfied ``near the top'' of
the portraits of the elements in $\G$ (and therefore in $\GG$ as
well). These constraints, if satisfied by an element $g$ in
$\Aut(\tree)$ at each of its sections (see below for details),
guarantee that $g$ belongs to the closure $\GG$. The constraints may
be viewed as an effective version of the more conceptual description
given by Grigorchuk in~\cite{grigorchuk:unsolved}.

Grigorchuk group can be viewed as a group of automorphisms of the
binary rooted tree $\tree$. The vertices of the rooted binary tree
$\tree$ are in bijective correspondence with the finite words over
$X=\{0,1\}$. The empty word $\emptyset$ is the root, the set $X^n$
of words of length $n$ over $X$ constitutes level $n$ in the tree
and every vertex $u$ at level $n$ has two children at level $n+1$,
namely $u0$ and $u1$. The group $\Aut(\tree)$ of automorphisms of
$\tree$ decomposes algebraically as
\begin{equation}\label{decomposition}
 \Aut(\tree) = (\Aut(\tree) \times \Aut(\tree)) \rtimes \Sym(2),
\end{equation}
where $\Sym(2) = \{1,\sigma\}= \{(),(01)\}$ is the symmetric group
of order 2 acting on $\Aut(\tree) \times \Aut(\tree)$ by permuting
the coordinates. The normal subgroup $\Aut(\tree) \times
\Aut(\tree)$ is the stabilizer of the first level of $\tree$ and the
elements in $\Aut(\tree)$ of the form $g=(g_0,g_1)$ act on $\tree$
by
\[ (0w)^g = 0w^{g_0} , \qquad (1w)^g = 1w^{g_1}, \]
while the elements of the form $g = (g_0,g_1) \sigma$ act by
\[ (0w)^g = 1w^{g_0} , \qquad (1w)^g = 0w^{g_1},  \]
for any word $w$ over $X$. The automorphisms $g_0$ and $g_1$ in the
decomposition $g=(g_0,g_1)\sigma^\varepsilon$ of $g$, where
$\varepsilon$ is 0 or 1, are called \emph{sections} of $g$ at the
vertices 0 and 1, respectively. This definition is recursively
extended to a notion of a section of $g$ at any vertex of $\tree$ by
declaring $g_\emptyset=g$ and $g_{ux} = (g_u)_x$, for $u$ a word
over $X$ and $x$ a letter in $X$.

\emph{Grigorchuk group} $\G$ is, by definition, the group generated
by the automorphisms $a$, $b$, $c$ and $d$ of $\tree$, whose
decompositions, in the sense of \eqref{decomposition}, are given by
\begin{align}
 a &= (1,1)\sigma, \notag\\
 b &= (a,c), \notag\\
 c &= (a,d), \label{c}\\
 d &= (1,b). \notag
\end{align}
Therefore, the action of $a$, $b$, $c$, and $d$ on $\tree$ is given
by
\begin{alignat*}{2}
 (0w)^a &= 1w,  & \qquad (1w)^a &= 0w, \\
(0w)^b &= 0w^a, & \qquad (1w)^b &= 1w^c, \\
(0w)^c &= 0w^a, & \qquad (1w)^c &= 1w^d, \\
(0w)^d &= 0w ,  & \qquad (1w)^d &= 1w^b,
\end{alignat*}
for any word $w$ over $X$. It is easy to establish that
\[
 a^2=b^2=c^2=d^2=1, \qquad bc=cb=d, \qquad bd=db=c, \qquad cd=dc=b.
\]
These relations are called simple relations in $\G$.  The stabilizer
$\Stab_\G(X)$ in $\G$ of level 1 in $\tree$ is
\[ \Stab_\G(X) = \langle b,c,d,aba,aca,ada \rangle, \]
and the decompositions of $aba$, $aca$ and $ada$ are given by
\begin{align}
 aba &= (c,a), \notag\\
 aca &= (d,a), \label{aca} \\
 ada &= (b,1). \notag
\end{align}

The decomposition formulae given in \eqref{c} and \eqref{aca} and
the simple relation $aa=1$ are sufficient to calculate the
decomposition of any element in $\G$. For example,
\[ abdabac = aba \ ada \ b \ aca \ a = (c,a)(b,1)(a,c)(d,a)(1,1)\sigma = (cbad,aca)\sigma.\]
Of course, we could make use of the other simple relations to write
either
\[ abdabac = \dots = (cbad,aca)\sigma = (dad,aca)\sigma.\]
or
\[ abdabac = acabac = aca \ b \ aca \ a  = (d,a)(a,c)(d,a)(1,1)\sigma = (dad,aca)\sigma, \]
but this will not be necessary for our purposes (and would, in fact,
be counterproductive in the proof of one of our lemmata).

Let $g$ be an arbitrary element in $\Aut(\tree)$. The
\emph{portrait} of $g$ is the binary rooted tree $\tree$ with
additional \emph{decoration} on the vertices defined recursively as
follows. If $g = (g_0,g_1)$ stabilizes level 1 in $\tree$ then the
portrait of $g$ consists of the portrait of $g_0$ hanging below the
vertex 0, the portrait of $g_1$ hanging below the vertex 1 and the
root, which is decorated by 0. If $g=(g_0,g_1)\sigma$ does not
stabilize level 1 (i.e., it is active at the root) the portrait
looks the same as in the previous case, except that the root is
decorated by 1. Thus, the portrait of $g$ is the binary tree $\tree$
with additional decoration $\alpha_u(g)$ on each vertex $u$, which
is equal to 0 or 1 depending on whether $g$ is active at the vertex
$u$ or not.

For every vertex $xy$ on level 2, define
\[
 \beta_{xy}(g) = \alpha_{xy}(g) + \alpha_{x\bar{y}0}(g) + \alpha_{x\bar{y}1}(g),
\]
where the addition is performed modulo 2 and $\bar{y}$ denotes the
letter in $\{0,1\}$ different from $y$. When $g$ is assumed, the
notation $\alpha_u(g)$ and $\beta_u(g)$ is simplified to $\alpha_u$
and $\beta_u$.

\begin{theorem}\label{t-constraints}
For any element $g$ in Grigorchuk group $\G$ the portrait decoration
satisfies the following constraints.

(i) If $\alpha_0=\alpha_1=0$ then
\[ \beta_{00} = \beta_{11} = \beta_{01} = \beta_{10}. \]

(ii) If $\alpha_0=0$ and $\alpha_1=1$ then
\[ \beta_{00} \neq \beta_{11} = \beta_{01} = \beta_{10}. \]

(iii) If $\alpha_0=1$ and $\alpha_1=0$ then
\[ \beta_{00} = \beta_{11} = \beta_{01} \neq \beta_{10}. \]

(iv) If $\alpha_0=\alpha_1=1$ then
\[ \beta_{00} = \beta_{11} \neq \beta_{01} = \beta_{10}. \]
\end{theorem}

We say that an automorphism $g$ in $\Aut(\tree)$ \emph{simulates}
$\G$ if its portrait decoration satisfies the constraints in Theorem
1. Recall that the pro-finite group $\Aut(\tree)$ is a metric space
with a natural metric defined as follows. The distance between two
tree automorphisms $g$ and $h$ is $\frac{1}{2^n}$, where $n$ is the
largest integer for which $g$ and $h$ agree on all words of length
$n$. In other words, if $g^{-1}h$ belongs to
$\Stab_{Aut(\tree)}(X^n)$, but not to
$\Stab_{\Aut(\tree)}(X^{n+1})$, then the distance between $g$ and
$h$ is $\frac{1}{2^n}$.

\begin{theorem}\label{closure-description}
Let $g$ be a binary tree automorphism. The following conditions are
equivalent.

(i) $g$ belongs to the closure $\GG$ of Grigorchuk group $\G$ in the
pro-finite group $\Aut(\tree)$.

(ii) all sections of $g$ simulate $\G$.

(iii) the distance from any section of $g$ to $\G$ in the metric
space $\Aut(\tree)$ is at most $\frac{1}{16}$.
\end{theorem}

The proofs will follow from some combinatorial observations on the
structure of words representing elements in $\G$. Before the proofs,
we consider some examples.

\begin{example}\label{strategy}
We show how Theorem~\ref{t-constraints} and
Theorem~\ref{closure-description} can be used to construct elements
in the closure $\GG$.

The constraints in Theorem~\ref{t-constraints} imply that exactly
$2^{12}$ different portrait decorations are possible on levels 0
through 3 for elements in $\G$. Indeed, assume that the portrait
decoration is already freely chosen on levels 0 through 2. In
particular, $\alpha_0$ and $\alpha_1$ are known. There are 8
vertices on level 3, but according to the constraints in
Theorem~\ref{t-constraints} we may choose the decoration freely only
on 5 of them. Namely, as soon as we chose the decoration for two
vertices with common parent, the values of $\beta_{00}$,
$\beta_{01}$, $\beta_{10}$ and $\beta_{11}$ are uniquely determined
and we may freely choose only the decoration on one of the vertices
in each of the 3 remaining pairs of vertices with common parent,
while the other is forced on us. For example, let us set
$\alpha_u=1$ for all $u$ on level 0 through level 2.
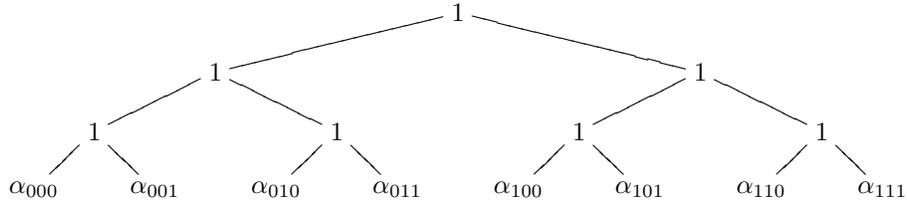
\[
 \xymatrix@C=5pt@R=10pt{
 &&&&&&& 1 \ar@{-}[dllll]\ar@{-}[drrrr]\\
 &&& 1 \ar@{-}[dll]\ar@{-}[drr] &&&&&&&& 1 \ar@{-}[dll]\ar@{-}[drr]\\
 & 1 \ar@{-}[dl]\ar@{-}[dr] &&&& 1 \ar@{-}[dl]\ar@{-}[dr] &&&&
   1 \ar@{-}[dl]\ar@{-}[dr] &&&& 1 \ar@{-}[dl]\ar@{-}[dr] \\
 \alpha_{000} && \alpha_{001} && \alpha_{010} && \alpha_{011} &&
 \alpha_{100} && \alpha_{101} && \alpha_{110} && \alpha_{111}
 }
\]
Further, for all vertices at level 3 whose label ends in 1 choose
$\alpha_{u1}=1$. Finally, choose $\alpha_{110}=1$. At this moment,
after making 12 free choices, we have that $\beta_{10}=1$ and
according to Theorem~\ref{t-constraints} we must have
$\beta_{01}=1$, $\beta_{00}=\beta_{11}=0$. In accordance with the
other choices already made on level 3, we must then have
\[
 \alpha_{000} = 1, \qquad \alpha_{010}=0, \qquad \alpha_{100}=0.
\]
We may now continue building a portrait of an element in $\GG$ by
extending (independently!) the left half and the right half of the
portrait one more level by following only the constraints imposed by
Theorem~\ref{t-constraints}.

As a general strategy (one can certainly choose a different one,
guided by any suitable purpose), we choose arbitrarily the
decoration on all vertices whose label ends in 1 or in 110, and then
fill in the decoration on the remaining vertices following
Theorem~\ref{t-constraints}.

If we continue our example above by decorating by 1 all vertices
whose label ends in 1 or in 110, we obtain
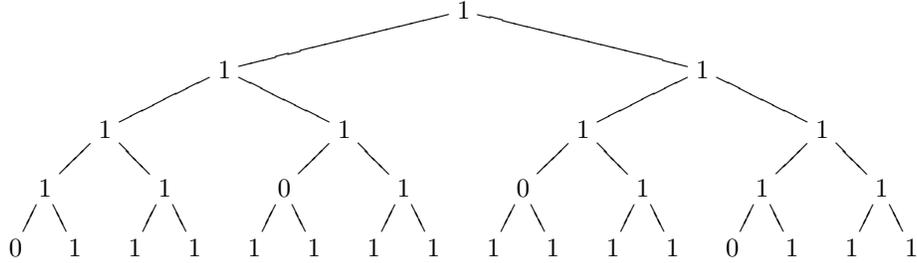
\[
 \xymatrix@C=0.3pt@R=10pt{
 &&&&&&&&&&&&&&& 1 \ar@{-}[dllllllll]\ar@{-}[drrrrrrrr]\\
 &&&&&&& 1 \ar@{-}[dllll]\ar@{-}[drrrr] &&&&&&&&&&&&&&&& 1 \ar@{-}[dllll]\ar@{-}[drrrr]\\
 &&& 1 \ar@{-}[dll]\ar@{-}[drr] &&&&&&&& 1 \ar@{-}[dll]\ar@{-}[drr] &&&&&&&&
     1 \ar@{-}[dll]\ar@{-}[drr] &&&&&&&& 1 \ar@{-}[dll]\ar@{-}[drr] \\
 & 1 \ar@{-}[dl]\ar@{-}[dr]&&&& 1 \ar@{-}[dl]\ar@{-}[dr]&&&& 0 \ar@{-}[dl]\ar@{-}[dr]&&&& 1 \ar@{-}[dl]\ar@{-}[dr]&&&&
   0 \ar@{-}[dl]\ar@{-}[dr]&&&& 1 \ar@{-}[dl]\ar@{-}[dr]&&&& 1 \ar@{-}[dl]\ar@{-}[dr]&&&& 1 \ar@{-}[dl]\ar@{-}[dr]\\
 0 && 1 && 1 && 1 && 1 && 1 && 1 && 1 && 1 && 1 && 1 && 1 && 0 && 1 && 1 && 1
 }
\]
Continuing in the same fashion (choosing 1 whenever possible) we
arrive at the portrait description of the element $f$ in $\GG$
defined by the following decomposition formulae
\begin{alignat*}{2}
 f    &= (\ell,r)\ &&\sigma, \\
 \ell &= (r,m)   \ &&\sigma, \\
 r    &= (m,r)   \ &&\sigma, \\
 m    &= (n,f)   \ &&\sigma, \\
 n    &= (r,m). \\
\end{alignat*}
There are many ways to see that $f$ (or any section of $f$) does not
belong to $\G$. Perhaps the easiest way is to observe that $f$ is
not a bounded automorphism, while all elements in $\G$ are bounded
automorphisms of $\tree$. Recall that, by definition, an
automorphism $g$ of $\tree$ is \emph{bounded} if the sum $\sum_{u
\in X^n} \alpha_u(g)$ is uniformly bounded, for all $n$. Note that
$f$ is defined by a 5-state automaton, which is not bounded, while
the automaton defining $\G$ is bounded. For more on groups of
automorphisms generated by automata
see~\cite{grigorchuk-n-s:automata}, and for bounded automorphisms
and bounded automata see~\cite{sidki:pol}.
\end{example}

To aid construction of elements in $\GG$ we provide the following
table of trees indicating the 8 possibilities for the values of
$\alpha_0$, $\alpha_1$, $\beta_{00}$, $\beta_{01}$, $\beta_{10}$,
and $\beta_{11}$.
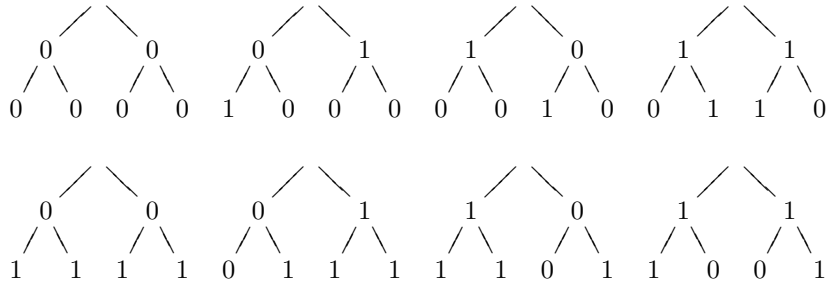
\[
 \xymatrix@C=0.3pt@R=10pt{
 &&&  \ar@{-}[dll]\ar@{-}[drr] &&&&&&&&  \ar@{-}[dll]\ar@{-}[drr] &&&&&&&&
      \ar@{-}[dll]\ar@{-}[drr] &&&&&&&&  \ar@{-}[dll]\ar@{-}[drr] \\
 & 0 \ar@{-}[dl]\ar@{-}[dr]&&&& 0 \ar@{-}[dl]\ar@{-}[dr]&&&& 0 \ar@{-}[dl]\ar@{-}[dr]&&&& 1 \ar@{-}[dl]\ar@{-}[dr]&&&&
   1 \ar@{-}[dl]\ar@{-}[dr]&&&& 0 \ar@{-}[dl]\ar@{-}[dr]&&&& 1 \ar@{-}[dl]\ar@{-}[dr]&&&& 1 \ar@{-}[dl]\ar@{-}[dr]\\
 0 && 0 && 0 && 0 && 1 && 0 && 0 && 0 && 0 && 0 && 1 && 0 && 0 && 1 && 1 && 0\\
 &&&  \ar@{-}[dll]\ar@{-}[drr] &&&&&&&&  \ar@{-}[dll]\ar@{-}[drr] &&&&&&&&
      \ar@{-}[dll]\ar@{-}[drr] &&&&&&&&  \ar@{-}[dll]\ar@{-}[drr] \\
 & 0 \ar@{-}[dl]\ar@{-}[dr]&&&& 0 \ar@{-}[dl]\ar@{-}[dr]&&&& 0 \ar@{-}[dl]\ar@{-}[dr]&&&& 1 \ar@{-}[dl]\ar@{-}[dr]&&&&
   1 \ar@{-}[dl]\ar@{-}[dr]&&&& 0 \ar@{-}[dl]\ar@{-}[dr]&&&& 1 \ar@{-}[dl]\ar@{-}[dr]&&&& 1 \ar@{-}[dl]\ar@{-}[dr]\\
 1 && 1 && 1 && 1 && 0 && 1 && 1 && 1 && 1 && 1 && 0 && 1 && 1 && 0 && 0 && 1\\
 }
\]
In each tree, $\alpha_0$, $\alpha_1$, $\beta_{00}$, $\beta_{01}$,
$\beta_{10}$, and $\beta_{11}$ are indicated in their respective
positions. To use the table, choose values for $\alpha_0$,
$\alpha_1$ and any one of $\beta_{00}$, $\beta_{01}$, $\beta_{10}$,
and $\beta_{11}$. The unique tree in the above table that agrees
with the chosen values provides the unique values for the remaining
3 parameters among $\beta_{00}$, $\beta_{01}$, $\beta_{10}$, and
$\beta_{11}$. For example, if $\alpha_0=\alpha_1=1$ and
$\beta_{11}=0$, the correct pattern in the table is the one in the
right upper corner, indicating that $\beta_{00}=0$, $\beta_{01}=1$
and $\beta_{10}=1$.

The following example provides additional ways to build elements in
$\GG$, which does not rely on Theorem~\ref{t-constraints} and
Theorem~\ref{closure-description}, but, rather, on the branch
structure of $\G$ (see~\cite{grigorchuk:jibg} for more details).
This approach does not produce all elements in $\GG$, but does
produce some that are easy to describe.

\begin{example}
Infinitely many elements in $\GG$ that are not in $\G$ are contained
in the following isomorphic copy of $K = [a,b]^\G \leq \G$ (recall
that $\G$ is a regular branch group over the normal closure $K$ of
$[a,b]$ in $\G$; see~\cite{grigorchuk;jibg} or~\cite{harpe:ggt} for
details). For each element $k \in K$ define an element $\bar{k}$ in
$\GG$ by
\[ \bar{k} = (k,\bar{k}). \]
The group $\bar{K} = \{\bar{k} \mid k \in K\}$ is canonically
isomorphic to $K$, but the intersection $\G \cap \bar{K}$ is
trivial.

More generally, an easy way to construct some (certainly not all)
elements in $\GG$ is to choose an infinite set of independent
vertices $V$ in $\tree$ (no vertex is  below some other vertex) and
associate to each such vertex an arbitrary element of $K$. The
automorphism of $\tree$ that is inactive at every vertex that does
not have a prefix in $V$ and whose sections at the vertices of $V$
are the assigned elements from $K$ is an element of $\GG$.
\end{example}

We now prove Theorem~\ref{t-constraints} and
Theorem~\ref{closure-description}.

Let $W$ be a word over $\{a,b,c,d\}$. The letters in $B=\{b,c,d\}$
are called $B$-letters. For a subset $C$ of $B$ the letters in $C$
are called $C$-letters. Denote by $N_C(W)$ the number of $C$-letters
occurring in $W$. An occurrence of a $B$-letter $\ell$ is called
even or odd depending on whether an even or odd number of $a$'s
appear before $\ell$ in $W$. For a parity $p \in \{0,1\}$ and a
subset $C$ of $B$, denote by $N^p_C(W)$ the number of $C$-letters of
parity $p$ in $W$. For parities $p,q \in \{0,1\}$, denote by
$N^{p,q}_{b,c}(W)$ the number of $\{b,c\}$-letters $\ell$ of parity
$p$ in $W$ such that the number of $\{b,c\}$-letters of parity
$\bar{p}$ that appear before $\ell$ in $W$ has parity $q$ (here
parity $\bar{p}$ denotes the parity different from $p$). For
example,
\begin{alignat*}{2}
 &N^1_{b,c}(a\mathbf{bc}aadabdbcad\mathbf{cb}d\mathbf{b}abdbc) &&= 5, \\
 &N^{1,1}_{b,c}(abcaadabdbcad\mathbf{cb}d\mathbf{b}abdbc) &&= 3, \\
 &N^{1,0}_{b,c}(a\mathbf{bc}aadabdbcadcbdbabdbc) &&= 2,
\end{alignat*}
where, in all examples, the letters that are counted are indicated
in boldface. When $W$ is assumed, the notation $N^p_C(W)$ and
$N^{p,q}_{b,c}(W)$ is simplified to $N^p_C$ and $N^{p,q}_{b,c}$.

\begin{lemma}\label{whoswho}
For any word $W$ over $\{a,b,c,d\}$ representing an element $g$ in
$\G$
\[
 N^{1,0}_{b,c} = \beta_{00}, \qquad N^{0,1}_{b,c} = \beta_{11}, \qquad
 N^{1,1}_{b,c} = \beta_{01}, \qquad N^{0,0}_{b,c} = \beta_{10},
\]
where all equalities are taken modulo 2.
\end{lemma}

\begin{proof}
Let the words $W_u$, for $u$ a word over $X$ of length at most 3,
represent the sections of $g$ at the corresponding vertices, and let
these words be obtained by decomposition from $W$, without applying
any simple relations other than $aa=1$ (i.e., no relations involving
$B$-letters are applied).

We have (modulo 2)
\begin{align*}
 \beta_{00} &= \alpha_{00} + \alpha_{010} + \alpha_{011} =
 \alpha_{00} + N_{b,c}(W_{01}) = \\
 &= N^0_{b,c}(W_0) + N^0_{b,d}(W_0) = N^0_{c,d}(W_0) = N^{1,0}_{b,c}(W).
\end{align*}

The other equalities are obtained in an analogous fashion.
\end{proof}

\begin{lemma}\label{l-constraints}
Let $W$ be a word over $\{a,b,c,d\}$. Modulo 2 we have

(i) if $N^0_{b,c} = 0$, then $N^{1,1}_{b,c} = N^{0,1}_{b,c} =
N^{0,0}_{b,c}$.

(ii) if $N^0_{b,c} = 1$, then $N^{1,0}_{b,c} = N^{0,1}_{b,c} \neq
N^{0,0}_{b,c}$.

(iii) if $N^1_{b,c} = 0$, then $N^{0,1}_{b,c} = N^{1,1}_{b,c} =
N^{1,0}_{b,c}$.

(iv) if $N^1_{b,c} = 1$, then $N^{0,0}_{b,c} = N^{1,1}_{b,c} \neq
N^{1,0}_{b,c}$.

\end{lemma}

\begin{proof}
(i) Assume $N^0_{b,c}$ is even.

The structure of the word $W$ may be represented schematically by
\[
 W = n_1 \ \ell_1 \ n_1' \ \ell_2 \ n_2 \ \ell_3 \ n_2' \ \dots \
     n_k \ \ell_{2k-1} \ n_k' \ \ell_{2k} \ n_{k+1},
\]
where $\ell_i$, $i=1,\dots,2k$ represent all the even occurrences of
$\{b,c\}$-letters in $W$ and the numbers $n_i,n_i'$ represent the
number of odd occurrences of $\{b,c\}$-letters between the
consecutive even occurrences of $\{b,c\}$-letters.

Then (modulo 2)
\[
 N^{1,1}_{b,c} = \sum_{i=1}^k n_i' =
 |\{ i \mid 1 \leq i \leq k, \ n_i' \text{ is odd }\}| = N^{0,1}_{b,c}.
\]
Indeed, for $i=1,\dots,k$, the $n_i'$ odd occurrences of
$\{b,c\}$-letters between $\ell_{2i-1}$ and $\ell_{2i}$ are preceded
by an odd number (exactly $2i-1$) of even occurrences of
$\{b,c\}$-letters. Thus $N^{1,1}_{b,c} = \sum_{i=1}^k n_i'$. On the
other hand, for $i=1,\dots,k$, whenever $n_i'$ is odd exactly one of
$\ell_{2i-1}$ and $\ell_{2i}$ is preceded by an odd number of odd
occurrences of $\{b,c\}$-letters, while whenever $n_i'$ is even
either both or none of $\ell_{2i-1}$ and $\ell_{2i}$ are preceded by
an odd number of odd occurrences of $\{b,c\}$-letters. Thus $|\{ i
\mid 1 \leq i \leq k, \ n_i' \text{ is odd }\}| = N^{0,1}_{b,c}$
modulo 2.

Since $N^{0,0}_{b,c}+N^{0,1}_{b,c} = N^0_{b,c}$ is even, we clearly
have $N^{0,0}_{b,c}=N^{0,1}_{b,c}$, modulo 2.

(ii) Assume $N^0_{b,c}$ is odd.

The structure of the word $W$ may be represented schematically by
\[
 W = n_1 \ \ell_1 \ n_1' \ \ell_2 \ n_2 \ \ell_3 \ n_2' \ \dots \
     n_{k-1}' \ \ell_{2k-2} \ n_k \ \ell_{2k-1} \ n_k' .
\]
where $\ell_i$, $i=1,\dots,2k-1$ represent all the even occurrences
of $\{b,c\}$-letters in $W$ and the numbers $n_i,n_i'$ represent the
number of odd occurrences of $\{b,c\}$-letters between the
consecutive even occurrences of $\{b,c\}$-letters.

Then (modulo 2)
\[
 N^{1,0}_{b,c} = \sum_{i=1}^k n_i =
 |\{ i \mid 1 \leq i \leq k, \ n_i \text{ is odd }\}| = N^{0,1}_{b,c}.
\]
Indeed, for $i=2,\dots,k$, the $n_i$ odd occurrences of
$\{b,c\}$-letters between $\ell_{2i-2}$ and $\ell_{2i-1}$ are
preceded by an even number (exactly $2i-2$) of even occurrences of
$\{b,c\}$-letters. In addition, the $n_1$ odd occurrences of
$\{b,c\}$-letters from the beginning of $W$ are preceded by no even
occurrences of  $\{b,c\}$-letters. Thus $N^{1,0}_{b,c} =
\sum_{i=1}^k n_i$. On the other hand, for $i=2,\dots,k$, whenever
$n_i$ is odd exactly one of $\ell_{2i-2}$ and $\ell_{2i-1}$ is
preceded by an odd number of odd occurrences of $\{b,c\}$-letters,
while whenever $n_i$ is even either both or none of $\ell_{2i-2}$
and $\ell_{2i-1}$ are preceded by an odd number of odd occurrences
of $\{b,c\}$-letters. In addition, whether $\ell_1$ is preceded by
an even or odd number of odd occurrences of $\{b,c\}$-letters
depends on the parity of $n_1$. Thus $|\{ i \mid 1 \leq i \leq k, \
n_i \text{ is odd }\}| = N^{0,1}_{b,c}$ modulo 2.

Since $N^{0,0}_{b,c}+N^{0,1}_{b,c} = N^0_{b,c}$ is odd, we clearly
have $N^{0,0}_{b,c} \neq N^{0,1}_{b,c}$, modulo 2.

(iii) and (iv) are analogous to (i) and (ii).
\end{proof}

\begin{proof}[Proof of Theorem~\ref{t-constraints}]
Follows directly from Lemma~\ref{whoswho},
Lemma~\ref{l-constraints}, and the observations $\alpha_0 =
N^0_{b,c}$ and $\alpha_1 = N^1_{b,c}$ modulo 2.
\end{proof}

\begin{proof}[Proof of Theorem~\ref{closure-description}]
We use the following (modification of the) description of the
elements in $\GG$ provided in~\cite{grigorchuk:unsolved}. A binary
tree automorphism $g$ belongs to $\GG$ if and only, for each section
$g_u$ of $g$, the portrait of $g_u$ agrees with the portrait of some
element in $\G$ up to and including level 3.

(i) is equivalent to (iii). Portraits of two automorphisms agree at
least up to level 3 if and only if their actions on the tree agree
at least up to level 4, which, in turn, is equivalent to the
condition that the distance between the two automorphisms is at most
$\frac{1}{2^4}=\frac{1}{16}$.

(i) implies (ii). If $g$ is in $\GG$, then the portrait of each
section $g_u$ of $g$ agrees with the portrait of some element in
$\G$ up to and including level 3. The portrait decorations of the
elements in $\G$ must satisfy the constraints in Theorem 1, and
therefore each section $g_u$ simulates $\G$.

(ii) implies (i). It is known that $\G/\Stab_\G(X^4)=2^{12}$. Thus,
for elements in $\G$, there are exactly $2^{12}$ possible portrait
decorations on level 0 through 3 . The constraints of
Theorem~\ref{t-constraints} provide for exactly $2^{12}$ different
decorations of the appropriate size (see the discussion in
Example~\ref{strategy}. Thus if a tree automorphism simulates $\G$,
then its portrait agrees with the portrait of an actual element in
$\G$ up to and including level 3.
\end{proof}

As another application, we offer a proof of the following result,
obtained by Grigorchuk in~\cite{grigorchuk:jibg}.

\begin{theorem}
The Hausdorff dimension of $\GG$ in $\Aut(\tree)$ is $\frac{5}{8}$.
\end{theorem}
\begin{proof}
It is known that the Hausdorff dimension can be calculated as the
limit
\[ \lim_{n \to \infty}
\frac{\log|\Stab_\GG(X^n)|}{\log|\Stab_{\Aut(\tree)}(X^n)|}, \]
comparing the sizes of the level stabilizers of $\GG$ and
$\Aut(\tree)$ (see~\cite{barnea-s:hausdorff}). Applying the strategy
of construction of elements in $\GG$ indicated in
Example~\ref{strategy}, it follows that, in the portrait of an
element $g$ in $\GG$, 5 out of 8 vertices at level 3 and below can
have any decoration we choose (0 or 1) and the other three have
uniquely determined decoration. Thus, the limit determining the
Hausdorff dimension of $\GG$ is $\frac{5}{8}$.
\end{proof}

\def\cprime{$'$}


\begin{thebibliography}{GNS00}

\bibitem[BS97]{barnea-s:hausdorff}
Yiftach Barnea and Aner Shalev.
\newblock Hausdorff dimension, pro-{$p$} groups, and {K}ac-{M}oody algebras.
\newblock {\em Trans. Amer. Math. Soc.}, 349(12):5073--5091, 1997.

\bibitem[dlH00]{harpe:ggt}
Pierre de~la Harpe.
\newblock {\em Topics in geometric group theory}.
\newblock Chicago Lectures in Mathematics. University of Chicago Press,
  Chicago, IL, 2000.

\bibitem[GNS00]{grigorchuk-n-s:automata}
R.~I. Grigorchuk, V.~V. Nekrashevich, and V.~I. Sushchanski{\u\i}.
\newblock Automata, dynamical systems, and groups.
\newblock {\em Tr. Mat. Inst. Steklova}, 231(Din. Sist., Avtom. i Beskon.
  Gruppy):134--214, 2000.

\bibitem[Gri80]{grigorchuk:burnside}
R.~I. Grigorchuk.
\newblock On {B}urnside's problem on periodic groups.
\newblock {\em Funktsional. Anal. i Prilozhen.}, 14(1):53--54, 1980.

\bibitem[Gri00]{grigorchuk:jibg}
R.~I. Grigorchuk.
\newblock Just infinite branch groups.
\newblock In Markus P. F. du~Sautoy Dan~Segal and Aner Shalev, editors, {\em
  New horizons in pro-$p$ groups}, pages 121--179. Birkh\"auser Boston, Boston,
  MA, 2000.

\bibitem[Gri05]{grigorchuk:unsolved}
R.~I. Grigorchuk.
\newblock Solved and unsolved problems around one group.
\newblock In {\em Progress in Mathematics}, volume 217, pages 117--217.
  Birkh\"auser, Basel, 2005.

\bibitem[Sid04]{sidki:pol}
Said Sidki.
\newblock Finite automata of polynomial growth do not generate a free group.
\newblock {\em Geom. Dedicata}, 108:193--204, 2004.

\end{thebibliography}
\end{document}